\documentclass{amsart}
\usepackage{graphicx}
\usepackage{amssymb}
\usepackage{amsfonts}
\swapnumbers
\newtheorem{theorem}{Theorem}[section]

\theoremstyle{definition}

\newtheorem{assumption}[theorem]{Assumption}
\newtheorem{remark}[theorem]{Remark}
\newtheorem*{acknowledgment}{Acknowledgment}

\numberwithin{equation}{section}
 \theoremstyle{plain}    
 
 \numberwithin{equation}{section} 
 \numberwithin{figure}{section} 
 \theoremstyle{plain}    
 \theoremstyle{plain}    
 \theoremstyle{remark}    
 \newtheorem*{acknowledgement*}{Acknowledgement} 

\newcommand{\cF}{{\mathcal F}}
\newcommand{\cG}{{\mathcal G}}

\newcommand{\cS}{{\mathcal S}}
\newcommand{\cT}{{\mathcal T}}

\newcommand{\te}{{\theta}}

\newcommand{\Om}{{\Omega}}
\newcommand{\om}{{\omega}}
\newcommand{\ve}{{\varepsilon}}
\newcommand{\del}{{\delta}}
\newcommand{\Del}{{\Delta}}
\newcommand{\gam}{{\gamma}}

\newcommand{\vf}{{\varphi}}
\newcommand{\vr}{{\varrho}}

\newcommand{\sig}{{\sigma}}
\newcommand{\al}{{\alpha}}
\newcommand{\be}{{\beta}}
\newcommand{\ka}{{\kappa}}
\newcommand{\la}{{\lambda}}

\newcommand{\vp}{{\varpi}}

\newcommand{\bbN}{{\mathbb N}}

\newcommand{\bbR}{{\mathbb R}}
\newcommand{\bbT}{{\mathbb T}}

\newcommand{\SRB}{{\mu^{\mbox{\tiny{SRB}}}}}

\begin{document}
\title[]{Nonconventional limit theorems in averaging}%
 \vskip 0.1cm 
 \author{ Yuri Kifer\\
\vskip 0.1cm
 Institute  of Mathematics\\
Hebrew University\\
Jerusalem, Israel}%
\address{
Institute of Mathematics, The Hebrew University, Jerusalem 91904, Israel}
\email{ kifer@math.huji.ac.il}%

\thanks{ }
\dedicatory{  }
 \date{\today}
\begin{abstract}\noindent
We consider "nonconventional" averaging setup in the form 
$\frac {dX^\epsilon(t)}{dt}=\epsilon B\big(X^\epsilon(t),\xi(q_1(t)),
\xi(q_2(t)),...,\xi(q_\ell(t))\big)$ where $\xi(t),t\geq 0$ is either a 
stochastic process or a dynamical system (i.e. then $\xi(t)=F^tx$) with
sufficiently fast mixing while 
$q_j(t)=\al_jt,\,\al_1<\al_2<...<\al_k$ and $q_j,\, j=k+1,...,\ell$ grow 
faster than linearly.
We show that the properly normalized error term in the "nonconventional"
averaging principle is asymptotically Gaussian. 
\end{abstract}
\maketitle
\markboth{Yu.Kifer}{Nonconventional averaging} 
\renewcommand{\theequation}{\arabic{section}.\arabic{equation}}
\pagenumbering{arabic}

\section{Introduction}\label{sec1}\setcounter{equation}{0}

Nonconventional ergodic theorems (see \cite{Fu}) known also after \cite{Be}
as polynomial ergodic theorems studied the limits of expressions having the 
form $1/N\sum_{n=1}^NF^{q_1(n)}f_1\cdots F^{q_\ell (n)}f_\ell$ where $F$ is a
weakly mixing measure preserving transformation, $f_i$'s are bounded measurable
functions and $q_i$'s are polynomials taking on integer values on the integers.
Originally, these results were motivated by applications to multiple recurrence
for dynamical systems taking functions $f_i$ being indicators of some measurable
sets and only convergence in the $L^2$-sense was dealt with but later \cite{As}
provided also almost sure convergence under additional conditions. Recently 
such results were extended in \cite{BLM} to the continuous time dynamical 
systems, i.e. to expressions of the form
\[
\frac 1{\cT}\int_0^{\cT}F^{q_1(t)}f_1\cdots F^{q_\ell (t)}f_\ell dt
\]
where $F^s$ is now an ergodic measure preserving flow.

In this paper we consider the averaging setup 
\begin{equation}\label{1.1}
X^\ve(n+1)=X^\ve(n)+\ve B(X^\ve(n),\xi(q_1(n)),...,\xi(q_\ell(n)))
\end{equation}
in the discrete time case and
\begin{equation}\label{1.2}
\frac {dX^\ve(t)}{dt}=\ve B(X^\ve(t),\xi(q_1(t)),...,\xi(q_\ell(t)))
\end{equation}
in the continuous time case with $\xi$ being either a stochastic process or 
having the form $\xi(s)=F^sf$ where $F^s$ is a dynamical system and $f$ is a 
function. Positive functions $q_1,...,q_\ell$ will satisfy certain conditions 
which will be specified in the next section, in particular, first $k$ of them
are linear while others grow faster than preceeding ones. An example where 
(\ref{1.2}) emerges is obtained when we consider a time dependent small 
perturbation of the oscillator equation
\begin{equation}\label{1.3}
\ddot x+\la^2x=\ve g(x,\dot x,t)
\end{equation}
where the force term $g$ depends on time in a random way $g(x,y,t)=g(x,y,
\xi(q_1(t)),...,\xi(q_\ell(t)))$. Then passing to the polar coordinates 
$(r,\phi)$ with $x=r\sin(\la(t-\phi))$ and $\dot x=\la r\cos(\la(t-\phi))$ 
the equation (\ref{1.3}) will be transformed into (\ref{1.2}). It seems 
reasonable that a random force may depend on versions of a same process or a 
dynamical system moving with different speeds which is what we have here.

As it is well known (see, for instance, \cite{SV}), if $B(x,y_1,...,y_\ell)$ 
is bounded and Lipschitz continuous in $x$ and the limit
\begin{equation}\label{1.4}
\bar B(x)=\lim_{\cT\to\infty}\frac 1\cT\int_0^\cT B(x,\xi(q_1(t)),...,
\xi(q_\ell(t)))dt
\end{equation}
exists then for any $S\geq 0$,
\begin{equation}\label{1.5}
\lim_{\ve\to 0}\sup_{0\leq t\leq S/\ve}|X^\ve(t)-\bar X^\ve(t)|=
\lim_{\ve\to 0}\sup_{0\leq t\leq S}|Z^\ve(t)-\bar Z(t)|=0
\end{equation}
where
\begin{equation}\label{1.6}
\frac {d\bar X^\ve(t)}{dt}=\ve\bar B(\bar X^\ve(t))\,\,\,\mbox{and}\,\,\,
Z^\ve(t)=X^\ve(t/\ve),\,\bar Z(t)=\bar X^\ve(t/\ve).
\end{equation}
In the discrete time case we have to take
\begin{equation}\label{1.7}
\bar B(x)=\lim_{N\to\infty}\frac 1N\sum_{n=0}^N B(x,\xi(q_1(n)),...,
\xi(q_\ell(n)))
\end{equation}
and (\ref{1.5}) remains true with $\bar X^\ve$ given by (\ref{1.6}) and
(\ref{1.7}). Almost everywhere limits in (\ref{1.4}) and (\ref{1.7}) can
 be obtained by nonconventional
pointwise ergodic theorems from \cite{BLM} and \cite{As}, respectively, in 
rather general circumstances in the dynamical systems case and under another 
set of conditions existence of such limits follows from \cite{Ki4}.

After nonconventional ergodic theorems (or in the probabilistic language laws 
of large numbers) are established the next natural step is to obtain central 
limit theorem type results which was accomplished in \cite{KV}. The averaging 
principle (\ref{1.5}) can be considered as an extension of the ergodic theorem
and the main goal of this paper is to extend also central limit theorem type 
results to the above nonconventional averaging setup in the spirit of what was
done in the standard (conventional) averaging case in \cite{Kh1} and \cite{Ki1}.
Central limit theorem type results turn in the averaging setup into assertions 
about Gaussian approximations of the slow motion $X^\ve$ given by (\ref{1.1}) 
or by (\ref{1.2}) where $\xi$ is a fast mixing stochastic process or a dynamical
system while unlike the standard (conventional) case we have the process $\xi$
taken simultaneously at different times $q_i(t)$ in the right hand side of
(\ref{1.1}) and (\ref{1.2}). 

We prove, first, our limit theorems for stochastic
processes under rather general conditions resembling the definition of 
mixingales (see \cite{ML1} and \cite{ML2}) and then check these
conditions for more familiar classes of stochastic processes and dynamical
systems. In \cite{KV} we imposed mixing assumptions in a standard way
relying on two parameter families of $\sig$-algebras (see \cite{Br}) while
our assumptions here use only filtrations (i.e. nondecreasing families)
of $\sig$-algebras which are easier to construct for various classes of
dynamical systems. As one of applications we check some form of 
our conditions for Anosov
flows which serve as fast motions in our nonconventional averaging setup
where we rely on the notion of Markov families from \cite{Do1} and \cite{Do2}.

At the end of the paper we discuss a fully coupled averaging setup in our 
nonconventional situation where already an averaging principle itself becomes a 
problem.

\begin{acknowledgment}
A part of this paper was done during my visit to the Fields Institute in
Toronto in June 2011 whose support and excellent working conditions I 
greatfully acknowledge.
\end{acknowledgment}

\section{Preliminaries and main results}\label{sec2}\setcounter{equation}{0}

Our setup consists of a  $\wp$-dimensional stochastic process
$\{\xi(t),\, t\geq 0\,\,\mbox{or}\,\, t=0,1,...\}$ on a probability space
$(\Om,\cF,Pr)$  together with a filtration of $\sig$-algebras 
$\cF_l\subset\cF,\, 0\leq l\leq\infty$ 
so that $\cF_l\subset\cF_{l'}$ if $l\leq l'$.
For convenience we extend the definitions of  $\cF_l$ given only for $l\geq 0$
to negative $l$ by defining  $\cF_l=\cF_0$ for $l<0$. In order to relax
required stationarity assumptions to some kind of weak "limiting stationarity"
our setup includes another probability measure $P$ on the space $(\Om,\cF)$.
Namely, we assume that the distribution of $\xi(t)$ with respect to $P$ does
not depend on $t$ and the joint distribution of
$\{\xi(t), \xi(t')\}$  for $t\geq t'$ depends only on $t-t'$ which can be
 written in the form 
\begin{equation}\label{2.1}
\xi(t)P=\mu\,\,\mbox{and}\,\,(\xi(t),\xi(t'))P=\mu_{t-t'}\,\,\mbox{for all}
\,\,t\geq t'
\end{equation}
where $\mu$ is a probability measure on $\bbR^\wp$ and $\mu_s,\, s\geq 0$ is 
a probability measure on $\bbR^\wp\times\bbR^\wp$.

Our setup relies on two probability measures $Pr$ and $P$ in order to include,
for instance, Markov processes $\xi(t)$ satisfying the Doeblin condition
(see \cite{IL} or \cite{Do}) starting at a fixed point or with another
noninvariant distribution. Then $Pr$ will be a corresponding probability
in the path space while $P$ will be the stationary probability constructed
by the initial distribution being the invariant measure of $\xi(t)$.
Usual mixing conditions for stochastic processes are formulated in terms of
a double parameter family of $\sig$-algebras via a dependence coefficient
between widely separated past and future $\sig$-algebras (cf. \cite{Br}
and \cite{KV}) but this approach
often is not convenient for applications to dynamical systems where natural
future $\sig$-algebras do not seem to exist unless an appropriate symbolic
representation is available. By this reason we formulate below a different 
set of mixing and approximation conditions for the process $\xi$ which seem
to be new and will enable us to treat some of dynamical systems models within a
class of stochastic processes satisfying our assumptions.

In order to avoid some of technicalities we restrict ourselves here mostly to
 bounded functions though our results can be obtained for more general
 classes of functions with polynomial growth supplemented by appropriate
 moment boundedness conditions similarly to \cite{KV}. For any function
 $g=g(\xi,\tilde\xi)$ on $\bbR^\wp\times\bbR^\wp$ introduce its H\" older norm 
\begin{equation}\label{2.2}
|g|_\ka=\sup\{ |g(\xi,\tilde\xi)|+\frac {|g(\xi,\tilde\xi)-g(\xi',\tilde\xi')|}
{|\xi-\xi'|^\ka+|\tilde\xi-\tilde\xi'|^\ka}:\, \xi\ne \xi',\, \xi\ne \xi'\}.
\end{equation}
Here and in what follows $|\psi-\tilde\psi|^\ka$ for two vectors $\psi=(\psi_1,
...,\psi_\vr)$ and $\tilde\psi=(\tilde\psi_1,...,\tilde\psi_\vr)$ denotes the
sum $\sum_{i=1}^\vr |\psi_i-\tilde\psi_i|^\ka$.
Next, for $p,q\geq 1$ and $s\geq 0$ we define a sort of a mixing coefficient
\begin{eqnarray}\label{2.3}
&\eta_{p,\ka,s}(n)=\sup_{t\geq 0}\big\{\big\| E\big(g(\xi(n+t),\xi(n+t+s))
|\cF_{[t]}\big)\\
&-E_Pg(\xi(n+t),\xi(n+t+s))\big\|_p:\, g=g(\xi,\tilde\xi),\, |g|_\ka\leq 1
\big\},\,\,\,\eta_{p,\ka}(n)=\eta_{p,\ka,0}(n)
\nonumber\end{eqnarray}
where $\|\cdot\|_p$ is the $L^p$-norm on the space $(\Om,\cF,Pr)$, $[\cdot]$ 
denotes the integral part and throughout this paper we write $E$ for the 
expectation with respect to $Pr$ and $E_P$ for the expectation with respect 
to $P$. We will need also an (one-sided) approximation coefficient
\begin{equation}\label{2.4}
\zeta_q(n)=\sup_{t\geq 0}\| E(\xi(t)|\cF_{[t]+n})-\xi(t)\|_q.
\end{equation}

\begin{assumption}\label{ass2.1}
Given $\ka\in(0,1]$ there exist $p,q\geq 1$ and $m,\del>0$ satisfying
\begin{equation}\label{2.5}
\gam_m=E|\xi(0)|^m<\infty,\, \frac 12\geq\frac 1p+\frac 2m+\frac {\del}q,\,
\del<\ka-\frac {\vr}p,\, \ka q>1
\end{equation}
with $\vr=(\ell-1)\wp$ and such that
\begin{equation}\label{2.6}
\sum_{n=0}^\infty n(\eta^{1-\frac {\vr}{p\te}}_{p,\ka}(n)+\zeta_q^\del(n))
<\infty\,\,\mbox{and}\,\,\lim_{n\to\infty}\eta_{p,\ka,s}(n)=0\,\,\mbox{for all}
\,\, s\geq 0,
\end{equation}
where $\frac {\vr}p<\te <\ka$.
\end{assumption}

Next, let $B=B(x,\xi)=(B^{(1)}(x,\xi),...,B^{(d)}(x,\xi)),$ $\xi=(\xi_1,...,
\xi_\ell)\in\bbR^{\ell\wp}$ be a 
$d$-vector function on $\bbR^d\times\bbR^{\ell\wp}$ such that for some
constant $K>0$ and all $x,\tilde x\in\bbR^d$, $\xi,\tilde\xi\in\bbR^{\ell\wp}$,
$i,j,l=1,...,d$,
\begin{eqnarray}\label{2.7}
&|B^{(i)}(x,\xi)|\leq K,\,|B^{(i)}(x,\xi)-B^{(i)}(\tilde x,\tilde\xi)|\leq
K(|x-\tilde x|+\sum_{j=1}^\ell|\xi_j-\tilde\xi_j|^\ka)\\
&\mbox{and}\,\,\quad\big\vert\frac {\partial B^{(i)}(x,\xi)}{\partial x_j}
\big\vert\leq K,\,\big\vert\frac {\partial^2 B^{(i)}(x,\xi)}{\partial x_j
\partial x_l}\big\vert\leq K.
\nonumber\end{eqnarray}
We will be interested in the central limit theorem type results as $\ve\to 0$
for the solution $X^\ve(t)=X^\ve_x(t)$ of the equation
\begin{equation}\label{2.8}
\frac {dX^\ve(t)}{dt}=\ve B\big( X^\ve(t),\xi(q_1(t)),\xi(q_2(t)),...,
\xi(q_\ell(t))\big),\, X^\ve_x(0)=x,\, t\in[0,\cT/\ve]
\end{equation}
where $q_1(t)<q_2(t)<\dots<q_\ell(t),\, t>0$ are increasing functions such
that $q_j(t)=\al_jt$ for $j\leq k<\ell$ with $\al_1<\al_2<\dots <\al_k$ 
whereas the remaining $q_j's$ grow faster in $t$. Namely, we assume 
similarly to \cite{KV} that for any $\gam>0$ and $k+1\leq i\leq\ell$,
\begin{equation}\label{2.9}
\lim_{t\to\infty}(q_i(t+\gam)-q_i(t))=\infty
\end{equation}
and
\begin{equation}\label{2.10}
\lim_{t\to\infty}(q_i(\gam t)-q_{i-1}(t))=\infty.
\end{equation}
 
Set 
\begin{equation}\label{2.11}
\bar B(x)=\int B(x,\xi_1,...,\xi_\ell)d\mu(\xi_1)\cdots d\mu(\xi_\ell).
\end{equation}
We consider also the solution $\bar X^\ve(t)=\bar X^\ve_x(t)$ of the averaged
 equation 
 \begin{equation}\label{2.12}
 \frac {d\bar X^\ve(t)}{dt}=\ve\bar B(\bar X^\ve(t)),\,\,\bar X^\ve_x(0)=x.
 \end{equation}
 It will be convenient to denote $Z^\ve(t)=X^\ve(t/\ve)$, $\bar Z(t)=
 \bar X^\ve(t/\ve)$ and to introduce $Y^\ve(t)=Y^\ve_y(t)$ by
 \begin{equation}\label{2.13}
Y_y^\ve(t)=y+\int_0^t B\big(\bar Z(s),\xi(q_1(s/\ve)),\xi(q_2(s/\ve)),...,
\xi(q_\ell(s/\ve))\big)ds.
\end{equation}
\begin{theorem}\label{thm2.2}
Suppose that (\ref{2.7}), (\ref{2.9}), (\ref{2.10}) and Assumption \ref{ass2.1}
hold true. Then the family of processes $G^\ve(t)=\ve^{-1/2}(Y^\ve_z(t) -
\bar Z_z(t)),\, t\in[0,\cT]$ converges weakly as $\ve\to 0$ to a Gaussian
process $G^0(t),\, t\in[0,\cT]$ having not necessarily independent increments
(see an example in \cite{KV}) with covariances of its components $G^0(t)=
(G^{0,1}(t),...,G^{0,d}(t))$ having the form $EG^{0,l}(s)G^{0,m}(t))=
\int_0^{\min(s,t)}A^{l,m}(u)du$ with the matrix function $\{A^{l,m}(u),\,
1\leq l,m\leq d\}$  computed in Section \ref{sec4}. Furthermore, the
family of processes $Q^\ve(t)=\ve^{-1/2}(Z^\ve(t)-\bar Z(t)),\, t\in[0,\cT]$
converges weakly as $\ve\to 0$ to a Gaussian process $Q^0(t),\, t\in[0,\cT]$
which solves the equation
\begin{equation}\label{2.14}
Q^0(t)=G^0(t)+\int_0^t\nabla\bar B(\bar Z(s))Q^0(s)ds.
\end{equation}

In the discrete time setup (\ref{1.1}) the similar results hold true assuming
that $q_i$'s take on integer values on integers, $\gam$ in (\ref{2.9}) is
replaced by 1, $\al_i$ is replaced by $i$ for $i=1,...,k$ and defining
$Z^\ve(t)=X^\ve([t/\ve])$ together with $Y^\ve=Y^\ve_y$ given by
\begin{equation}\label{2.15}
Y_y^\ve(t)=y+\int_0^t B\big(\bar Z(s),\xi(q_1([s/\ve])),\xi(q_2([s/\ve])),...,
\xi(q_\ell([s/\ve]))\big)ds
\end{equation}
while leaving all other definitions and assumptions the same as above.
\end{theorem}

Observe that we work with $\bar B$ defined by (\ref{2.11}) but in our
circumstances the central limit theorem type results from \cite{KV} imply
also (\ref{1.4}) and (\ref{1.7}), at least, in the $L^2$-sense while a
nonconventional law of large numbers from \cite{Ki4} and (under stationarity
assumptions) pointwise nonconventional ergodic theorems from \cite{As}
and \cite{BLM} yield (\ref{1.4}) and (\ref{1.7}) also for the almost
sure convergence. Note also that we need the full strength of (\ref{2.6})
only for one argument in Section \ref{sec4} borrowed from \cite{Kh1} but
for a standard limit theorem not in the averaging setup, i.e. when $B(x,\xi_1,
...,\xi_\ell)=B(\xi_1,...,\xi_\ell)$ does not depend on $x$, it suffices to
require only summability of the expression in brackets in (\ref{2.6}).

An important point in the proof of the first part of Theorem \ref{thm2.2}
is to introduce the representation
\begin{equation}\label{2.16}
B(x,\xi)=\bar B(x)+B_1(x,\xi_1)+\cdots +B_\ell(x,\xi_1,...,\xi_\ell)
\end{equation}
where $\xi=(\xi_1,...,\xi_\ell)$ and for $i<\ell$,
\begin{eqnarray}\label{2.17}
&B_i(x,\xi_1,...,\xi_\ell)=\\
&\int B(x,\xi_1,...,\xi_\ell)d\mu(\xi_{i+1})\cdots 
d\mu(\xi_\ell)-\int B(x,\xi_1,...,\xi_\ell)d\mu(\xi_{i})\cdots d\mu(\xi_\ell)
\nonumber\end{eqnarray}
while
\begin{equation}\label{2.18}
B_\ell(x,\xi_1,...,\xi_\ell)=B(x,\xi_1,...,\xi_\ell)-\int B(x,\xi_1,...,
\xi_\ell)d\mu(\xi_\ell).
\end{equation}
Next, we introduce
\begin{eqnarray}\label{2.19}
&\quad\,\,\,\, Y_i^\ve(t)=\int_0^{t/\al_i}B_i\big( \bar Z(s),\xi(q_1(s/\ve)),
\xi(q_2(s/\ve)),...,\xi(q_\ell(s/\ve))\big)ds,\,\, i=1,...,k,\\
&Y_i^\ve(t)=\int_0^tB_i\big( \bar Z(s),\xi(q_1(s/\ve)),\xi(q_2(s/\ve)),
...,\xi(q_\ell(s/\ve))\big)ds,\, i=k+1,...,\ell
\nonumber\end{eqnarray}
and $Y_0^\ve(t)=Y^\ve_{0,y}(t)=y+\int_0^tB_i(\bar Z(s))ds$. Thus $Y^\ve_y$
from (\ref{2.13}) has the representation
\begin{equation}\label{2.20}
Y_y^\ve(t)=Y^\ve_0(t)+\sum_{i=1}^kY_i^\ve(\al_it)+\sum_{i=k+1}^\ell Y^\ve_i(t).
\end{equation}
We consider also $X_0^\ve(t)=\bar X^\ve(t)$, $X^\ve_i(t)=X^\ve_{i,x}(t)=x+
\ve\int_0^tB_i\big(X^\ve_i(s),\xi(q_1(s)),...,\xi(q_\ell(s))\big)ds$
 and $Z^\ve_i(t)=X^\ve_i(t/\ve)$ for all $i\geq 0$.
For $i\geq 1$ set also
\begin{equation}\label{2.21}
G^\ve_i(t)=\ve^{-1/2}Y^\ve_i(t)\,\,\mbox{and}\,\, Q^\ve_i(t)=\ve^{-1/2}
Z^\ve_i(t).
\end{equation}

Relying on martingale approximations (which also can be done employing
 mixingales from \cite{ML1} and \cite{ML2}) we will show that any linear 
 combination $\sum_{i=1}^k\la_iG_i^\ve$ converges weakly as $\ve\to 0$ to
 a Gaussian process $\sum_{i=1}^k\la_iG_i^0$. It turns out that in the
  continuous time case each $G^\ve_i,\, i=k+1,...,\ell$ converges weakly
  as $\ve\to 0$ to zero, and so the processes $Y^\ve_i,\, i>k$ do not play
  any role in the limit. It follows that $G^\ve$ converges weakly to a
  Gaussian process $G^0$ such that $G(t)=\sum_{i=1}^k\la_iG_i^0(\al_it)$.
  On the other hand, in the discrete time case each $G^\ve_i,\, i>k$ cannot
  be disregarded, in general, and it converges weakly as $\ve\to 0$ to
  a Gaussian process $G^0_i$ which is independent of any other $G^0_j$. 
  The above difference between discrete and continuous time cases is due
  to the different natural forms of the assumption (\ref{2.9}) in these two
  cases. These arguments yield the first part of Theorem \ref{thm2.2}
 while its second part concerning convergence of $Q^\ve$ as $\ve\to 0$ is
  proved via some Taylor expansion and approximation arguments.
  
  In order to clarify the role of the coefficients $\eta_{p,\ka}$ and 
  $\zeta_q$ we compare them with the more familiar mixing and approximation 
  coefficients defined via a two parameter family of $\sig$-algebras
  $\cG_{s,t}\in\cF,\, -\infty\leq s\leq t\leq\infty$ by
 \begin{equation}\label{2.22}
 \vp_p(n)=\sup_{s\geq 0,g}\big\{\| E(g|\cG_{-\infty,s})-E_Pg\|_p:\, g\,\,
 \mbox{is}\,\,\cG_{s+n,\infty}-\mbox{measurable and}\,\, |g|\leq 1\big\}
 \end{equation}
 and
 \begin{equation}\label{2.23}
 \be_q(n)=\sup_{t\geq 0}\| E(\xi(t)|\cG_{t-n,t+n})-\xi(t)\|_q,
 \end{equation}
 respectively, where $\cG_{st}\subset\cG_{s't'}$ if $s'\leq s$ and $t'\geq t$.
 Then setting $\cF_l=\cG_{-\infty,l}$ we obtain by the contraction property 
 of conditional expectations that
 \begin{eqnarray}\label{2.24}
 &\be_q(n)\geq\sup_{t\geq 0}\| E(\xi(t)|\cG_{t-n,t+n})- \xi(t)+\xi(t)-
 E(\xi(t)|\cG_{-\infty,[t]+n+1})\|_q\\
 &\geq\zeta_q(n+1)-\be_q(n)\,\,\mbox{i.e.}\,\,\be_q(n)\geq\frac 12\zeta_q(n+1).
 \nonumber\end{eqnarray}
 Furthermore, 
 \begin{eqnarray*}
 &\big\| g(\xi(n+t),\xi(n+t+s))-g\big(E(\xi(n+t)|\cG_{n+t-[n/2],n+t+[n/2]}),\\
 &E(\xi(n+t+s)|\cG_{n+t+s-[n/2],n+t+s+[n/2]})\big)\big\|_p
 \leq 2|g|_\ka\be^\ka_{p\ka}([n/2]),
 \end{eqnarray*}
 and so
 \begin{equation}\label{2.25}
 \eta_{p,\ka}(n)\leq (\vp_p([n/2])+2\be_{p\ka}^\ka([n/2]))|g|_\ka.
 \end{equation}
 Thus, appropriate conditions on decay of coefficients $\vp_p$ and $\be_q$
 as in \cite{KV} yield corresponding conditions on $\eta_{p,\ka}$ and 
 $\zeta_q$. The other direction does not hold true but still it turns out 
 that most of the technique from \cite{KV} can be employed in our 
 circumstances, as well.
 
The conditions of Theorem \ref{thm2.2} hold true for many important stochastic
processes. In the continuous time case they are satisfied when, for instance,
$\xi(t)=f(\Xi(t))$ where $\Xi(t)$ is either an irreducible continuous time
finite state Markov chain or a nondegenerate diffusion process on a compact
manifold while $f$ is a H\" older continuous vector function. In the discrete
time case we can take, for instance, $\xi(n)=f(\Xi(n))$ with $\Xi(n)$ being
a Markov chain satisfying the Doeblin condition (see, for instance, \cite{IL},
p.p. 367--368). In all these examples $\eta_{p,\ka}(n)$ and $\zeta_q(n)$
decay in $n$ exponentially fast while (\ref{2.6}) requires much less. In
fact, in both cases $\xi(t)$ may depend on whole paths of a Markov process
$\Xi$ assuming only certain weak dependence on their tails.

Important classes of processes satisfying our conditions come from
dynamical systems. In Section \ref{sec6} we take $\xi(t)=\xi(t,z)
=g(F^tz)$ where $F^t$ is a $C^2$ Anosov flow (see \cite{KH}) on a 
compact manifold $M$ whose stable and unstable foliations are jointly
 nonintegrable and $g$ is a H\" older continuous $\wp$-vector function 
 on $M$. It turns out that if we take the initial point $z$ on an
 element $S$ of a Markov family (see Section \ref{sec6}) introduced in
 \cite{Do1} distributed there at random according to a probability
 measure equivalent to the volume on $S$ then Assumption \ref{ass2.1}
 can be verified. This does not yield though a desirable limit theorem
 where the initial point is taken at random on the whole manifold $M$
 distributed according to the Sinai-Ruelle-Bowen (SRB) measure (or the
 normalized Riemannian volume). 
We observe that a suspension representation of Anosov flows employed in
\cite{Ki1} to derive limit theorems in the conventional averaging setup
does not work in our situation because $F^{q_i(t)}x,\, i=1,...,\ell$ arrive
 at the ceiling of the suspension at different times for
  different $i$'s. 
  
  In the discrete time case there are several important classes of dynamical
  systems where our conditions can be verified. First,
  for transformations where symbolic representations via Markov
  partitions are available (Axiom A diffeomorphisms (see \cite{Bo}) and 
  expanding endomorphisms, some one-dimensional
  maps e.g. the Gauss map (see \cite{He}) etc.) we can rely on standard mixing
  and approximation assumptions based on two parameter families of 
  $\sig$-algebras as in (\ref{2.22}) and (\ref{2.23}). On the other hand,
  for many transformations Markov partitions are not available but still it
  is possible to construct one parameter increasing or decreasing filtration
  of $\sig$-algebras so that our conditions can be verified. For some classes
  of noninvertible transformations $F$ it is possible to choose an appropriate
  initial $\sig$-algebra $\cF_0$ such that $F^{-1}\cF_0\subset\cF_0$ and then
  to define a decreasing filtration $\cF_i=F^{-i}\cF_0$ (see \cite{Li}
  and \cite{FMT0}). Passing to the natural extension as in Remark 3.12 of 
  \cite{FMT0} we can turn to an increasing filtration and to verify our
  conditions.
  On the other hand, our results can be derived under appropriate conditions
  with respect to decreasing families of $\sig$-algebras. Namely, let
  $\cF\supset\cF_0\supset\cF_1\supset\cF_2\supset\cdots$ and define mixing
  and approximation coefficients by
  \begin{eqnarray}\label{2.26}
 &\eta_{p,\ka,s}(n)=\sup_{t\geq s}\big\{\big\| E\big(g(\xi(t),\xi(t-s))
 |\cF_{[t]+n}\big)\\
 &-E_Pg(\xi(t),\xi(t-s))\big\|_p:\, g=g(\xi,\tilde\xi),\, |g|_\ka\leq 1\big\},
 \,\,\,\eta_{p,\ka}(n)=\eta_{p,\ka,0}(n)
 \nonumber\end{eqnarray}
and 
\begin{equation}\label{2.27}
\zeta_q(n)=\sup_{t\geq n}\| E(\xi(t)|\cF_{[t]-n})-\xi(t)\|_q.
\end{equation}
Then under Assumption \ref{ass2.1} we can rely on estimates of Section 
\ref{sec3} below and in place of  martingales there arrive at reverse 
martingales and to use a limit theorem for the latter.

\begin{remark}\label{rem2.3}
  If $\bar B\equiv 0$ then according to Theorem \ref{thm2.2} the process
  $X^\ve(t)$ is very close to its initial point on the time interval of order
  $1/\ve$. Thus, in order to see fluctuations of order 1 it makes sense to
  consider longer time and to deal with $V^\ve(t)=X^\ve(t/\ve)$. Under the
  stronger condition $\int B(x,\xi_1,...,\xi_\ell)d\mu(\xi_\ell)\equiv 0$
  it is not difficult to mimic the proofs in \cite{Kh2} and \cite{Bor}
  relying on the technique of Sections \ref{sec3} and \ref{sec4} below
  in order to obtain that $V^\ve(t),\, t\in[0,T]$ converges weakly as 
  $\ve\to 0$ to a diffusion process with parameters obtained in the same way
  as in \cite{Kh2} and \cite{Bor}. It is not clear whether, in general, this
 result still holds true assuming only that $\bar B\equiv 0$. Though most of
 the required estimates still go through in the latter case a convergence 
 of $V^\ve$ to a Markov process seems to be problematic in a general 
 nonconventional averaging setup.
 \end{remark}

 \section{Estimates and martingale approximation}\label{sec3}
\setcounter{equation}{0}

The proof of Theorem \ref{thm2.2} will employ a modification of the machinery
 developed in \cite{KV}. First, we have to study the asymptotical behavior as
 $\ve\to 0$ of
\begin{equation}\label{3.1}
G^\ve_i(t)=\sqrt\ve\int_0^{\tau_i(t)/\ve}B_i\big(\bar Z(\ve s),\xi(q_1(s)),
...,\xi(q_i(s))\big)ds
\end{equation}
which is obtained from the definition (\ref{2.21}) by the change of variables
$s\to s/\ve$ and where $\tau_i(t)=t/\al_i$ for $i=1,...,k$ and $\tau_i(t)=t$
for $i=k+1,...,\ell$. Observe that if $\frac 1{N+1}\leq\ve\leq\frac 1N$ and 
$N\geq 1$ then by (\ref{2.7}),
\begin{equation}\label{3.2}
|G^\ve_i(t)-G^{1/N}_i(t)|\leq \frac {2Ktd}{\sqrt N},
\end{equation}
and so it suffices to study the asymptotical behavior of $G^{1/N}_i$ as
$N\to\infty$. Set
\begin{equation}\label{3.3}
I_{i,N}(n)=\int_n^{n+1}B_i\big(\bar Z(s/N),\xi(q_1(s)),...,\xi(q_i(s))\big)ds.
\end{equation}
In view of (\ref{2.7}) the asymptotical behavior of $G^{1/N}_i$ as 
$N\to\infty$ is the same as of $N^{-1/2}S_{i,N}(t)$ where
\begin{equation}\label{3.4}
S_{i,N}(t)=\sum_{n=0}^{[N\tau_i(t)]}I_{i,N}(n).
\end{equation}
There are two obstructions for applying directly the results of \cite{KV}
to the sum (\ref{3.4}). First, unlike \cite{KV} the integrand in (\ref{3.3})
depends on the "slow time" $s/N$. Secondly, our mixing and approximation
coefficients look differently from the corresponding coefficients in 
 \cite{KV}. Still, it turns out that these obstructions can
be dealt with and after minor modifications the method of \cite{KV} start
working in our situation, as well. Namely, the dependence on the "slow time"
being deterministic will not prevent us from making estimates similar to 
\cite{KV} while dependence of $I^N_i$ on $N$ will just require us to deal with
martingale arrays which creates no problems as long as we obtain appropriate
limits of variances and covariances. Concerning the second obstruction we
observe that one half of the approximation estimate from \cite{KV} is contained
in the coefficient $\zeta_p$ while another half is hidden in the coefficient
 $\eta_{p,\ka}$ which also suffices for required mixing estimates.
 
 We explain next more precisely why estimates similar to \cite{KV} hold true
 in our circumstances, as well. Let $f(\psi,\xi,\tilde\xi)$ be a function on 
 $\bbR^\vr\times\bbR^\wp\times\bbR^\wp$ such that for any $\psi,\psi'
 \in\bbR^\vr$ and $\xi,\tilde\xi,\xi',\tilde\xi'\in\bbR^\wp$,
 \begin{equation}\label{3.5}
 |f(\psi,\xi,\tilde\xi)-f(\psi',\xi',\tilde\xi')|\leq C(|\psi-\psi'|^\ka+
 |\xi-\xi'|^\ka+|y-y'|^\ka)\,\,\mbox{and}\,\, |f(\psi,\xi,\tilde\xi)|\leq C.
 \end{equation}
 Then setting $g(\psi)=E_Pf(\psi,\xi(0),\xi(s))$ we obtain from (\ref{2.1}) and 
 (\ref{2.3}) that for all $u,v\geq 0$ and $n\in\bbN$,
 \begin{equation}\label{3.6}
 \big\| E\big( f(\psi,\xi(n+u),\xi(n+u+v))|\cF_{[u]}\big)-g(\psi)\big\|_p\leq 
 C\eta_{p,\ka,v}(n).
 \end{equation}
 Let $h(\psi,\om)=E(f(\psi,\xi(n+u),\xi(n+u+v))|\cF_{[u]})-g(\psi)$. Then by 
 (\ref{3.5}) we can choose a version of $h(\psi,\om)$ such that with 
 probability  one simultaneously for all $\psi,\psi'\in\bbR^\vr$,
 \begin{equation}\label{3.7}
 |h(\psi,\om)-h(\psi',\om)|\leq 2C|\psi-\psi'|^\ka.
 \end{equation}
 Since, in addition, $\| h(\psi,\om)\|_p\leq C\eta_{p,\ka}(n)$ by (\ref{3.6})
 for all $\psi\in\bbR^\vr$, we obtain by Theorem 3.4 from \cite{KV} that for
 any random $\vr$-vector $\Psi=\Psi(\om)$,
 \begin{equation}\label{3.8}
 \| h(\Psi(\om),\om)\|_a\leq cC\big(\eta_{p,\ka,v}(n)\big)^{1-\frac {\vr}{p\te}}
 (1+\| \Psi\|_m)
 \end{equation}
 where $\frac {\vr}p<\te <\ka$, $\frac 1a\geq \frac 1p+\frac 1m$ and
 $c=c(\vr,p,\ka,\te)>0$ depends only on parameters in brackets. Since
 \begin{equation}\label{3.9}
 h(\tilde\Psi(\om),\om)=E\big( f(\tilde\Psi,\xi(n+u),\xi(n+u+v))|\cF_{[u]}\big)
 (\om) \,\,\mbox{a.s.}
 \end{equation}
 provided $\tilde \Psi$ is $\cF_{[u]}$-measurable we obtain from 
 (\ref{3.7})--(\ref{3.9}) together with the H\" older inequality  
 (cf. Corollary 3.6 in \cite{KV}) that,
 \begin{eqnarray}\label{3.10}
 &\big\| E\big( f(\Psi,\xi(n+u),\xi(n+u+v))|\cF_{[u]}\big)-g(\Psi)\big\|_a\\
 &\leq C(\eta_{p,\ka,v}(n))^{1-\frac {\vr}{p\te}}(1+\| \Psi\|_m)+2C\| \Psi-
 E(\Psi|\cF_{[u]})\|^\del_q\nonumber
 \end{eqnarray}
 provided $\frac 1a\geq \frac 1p+\frac 2m+\frac {\del}q$.
 
 We apply the above estimates in two cases. First, when $f(\psi,\xi,\tilde\xi)
 =f(\psi,\xi)=B_i(x,\xi_1,...,\xi_i)$ with $\psi=(\xi_1,...,\xi_{i-1})\in
 \bbR^{(i-1)\wp}$, $\xi=\xi_i\in\bbR^\wp$, $n=[(q_i(t)-q_{i-1}(t))/2]$,
 $u=q_i(t)-n$ and $\Psi=(\xi(q_1(t)),\xi(q_2(t)),...,\xi(q_{i-1}(t))$. In
 the second case $f(\psi,\xi,\tilde\xi)=B_i(x,\xi_1,...,\xi_i)B_j(y,\xi_1',
 ...,\xi_j')$ with $\psi=(\xi_1,...,\xi_{i-1},\xi_1',...,\xi_{j-1}')\in
 \bbR^{(i+j-2)\wp}$, $\xi=\xi_i$, $\tilde\xi=\xi_j'\in\bbR^\wp$, $n=
 \big[\big(\min(q_i(t),q_j(s))-\max(q_{i-1}(t),q_{j-1}(s))\big)/2\big]$ when 
 $n>0$, $u=\min(q_i(t),q_j(s))-n$ and $\Psi=\big(\xi(q_1(t)),...,
 \xi(q_{i-1}(t)),\xi(q_1(s)),...,\xi(q_{j-1}(s))\big)$. The estimates for the 
 first case are used for martingale approximations while the second case
 emerges when computing covariances.
 
 Since $\int B_i(x,\xi_1,...,\xi_{i-1},\xi_i)d\mu(\xi_i)=0$ we obtain by
 (\ref{3.10}) the estimate
 \begin{equation}\label{3.11}
 \big\| E\big(B_i(\xi(q_1(t)),...,\xi(q_i(t)))|\cF_{[q_i(t)]-n}\big)\big\|_a
 \leq C\big((\eta_{p,\ka}(n))^{1-\frac {\vr}{p\te}}+(\zeta_q(n))^\del\big)
 \end{equation}
 for some $C>0$ independent of $t$ where $n=n_i(t)=[(q_i(t)-q_{i-1}(t))/2]$. 
 Next, for any $x\in\bbR^\vr$, $\xi_1,...,\xi_{i-1}\in\bbR^\vr$ and $r=1,2,...$
 set
 \begin{eqnarray*}
 &B_{i,r}(x,\xi_1,...,\xi_{i-1},\xi(t))=E\big( B(x,\xi_1,...,\xi_{i-1},\xi(t))|
 \cF_{[t]+r}\big)\\
 &\mbox{and}\,\,\,\xi_r(t)=E\big(\xi(t)|\cF_{[t]+r}\big).
 \end{eqnarray*}
 Then by (\ref{2.4}) and (\ref{2.7}) together with the H\" older inequality,
 \begin{eqnarray}\label{3.12}
 &\big\| B_i(x,\xi_1,...,\xi_{i-1},\xi(t))-B_{i,r}(x,\xi_1,...,\xi_{i-1},
 \xi(t))\big\|_q\\
 &\leq 2\big\| B_i(x,\xi_1,...,\xi_{i-1},\xi(t))-B_i(x,\xi_1,...,\xi_{i-1},
 E(\xi(t)|\cF_{[t]+r}))\big\|_q\nonumber\\
 &\leq 2Kd\big\| |\xi(t)-E(\xi(t))|\cF_{[t]+r})|^\ka\big\|_q\leq 2Kd
 \zeta^\del_q(r).
 \nonumber\end{eqnarray}
 Moreover, similarly to Lemma 3.12 in \cite{KV} we obtain that
 \begin{equation}\label{3.13}
 \big\| B_i(x,\xi(q_1(t)),...,\xi(q_i(t)))-B_{i,r}(x,\xi_r(q_1(t)),...,
 \xi_r(q_{i-1}(t)))\big\|_a\leq c\zeta^\del_q(r)
 \end{equation}
 provided $\frac 1a\geq\frac 1p+\frac 2m+\frac {\del}q$ and $\del<\min(\ka,
 1-\frac d{p\ka})$ where $c=c(\del,a,p,q)>0$ depends only on the parameters
 in brackets. Set
 \[
 b^{l,m}_{ij}(x,y;s,t)=E\big( B^{(l)}_i(x,\xi(q_1(s)),...,\xi(q_i(s)))B^{(m)}_j
 (y,\xi(q_1(t)),...,\xi(q_j(t)))\big)
 \]
 where, recall, $B_i^{(l)}$ is the $l$-th component of the $d$-vector $B_i$.
 Now, by (\ref{2.7}), (\ref{3.11}) and (\ref{3.13}),
 \begin{equation}\label{3.14}
 |b^{l,m}_{ij}(x,y;s,t)|\leq C(\big((\eta_{p,\ka}(n))^{1-\frac {\vr}{p\te}}
 +(\zeta_q(n))^\del\big)
 \end{equation}
 where $C>0$ does not depend on $s,t\geq 0$ and $n=n_{ij}(s,t)=\max(\hat
 n_{ij}(s,t),\,\hat n_{ji}(t,s))$ with $\hat n_{ij}(s,t)=
 [\frac 12\min(q_i(s)-q_j(t),\, q_i(s)-q_{i-1}(s))]$.
 
 Now, set
 \begin{eqnarray}\label{3.15}
 &I_{i,N,r}(n)=\int_{n-1}^nB_{i,r}\big(\bar Z(s/N),\xi_r(q_1(s)),...,
 \xi_r(q_{i-1}(s))\big)ds,\, S_{i,N,r}(t)=\\
 &\sum_{n=1}^{[N\tau_i(t)]}I_{i,N,r}(n),\, R_{i,r}(m)=
 \sum_{l=m+1}^\infty E(I_{i,N,r}(l)|\cF_{m+r}),\, D_{i,N,r}(m)=\nonumber\\
 &I_{i,N,r}(m)+R_{i,r}(m)-R_{i,r}(m-1)\,\mbox{and}\,
 M_{i,N,r}(t)=\sum_{n=1}^{[N\tau_i(t)]}D_{i,N,r}(n).
 \nonumber\end{eqnarray}
 In view of (\ref{3.11}) applied with $a=2$ we see that the series for
 $R_{i,r}(m)$ converges in $L^2$, $D_{i,N,r}(m)$ is $\cF_{m+r}$-measurable
  and since $E(D_{i,N,r}(m)|\cF_{m-1+r})=0$
 we obtain that $\{ D_{i,N,r}(m),\cF_{m+r}\}_{0\leq m\leq [N\tau_i(T)]}$ is a 
 martingale differences array. Next, we proceed similarly to Sections 5 
 and 7 of \cite{KV} observing that the limiting behaviour of
 $N^{-1/2}S_{i,N,r}(t)$ as $N\to\infty$ is the same as of $N^{-1/2}
 M_{i,N,r}(t)$, then dealing with the latter by means of martingale limit
 theorems and, finally,  employing the representation
 \begin{equation}\label{3.16}
 S_{i,N}(t)=S_{i,N,1}(t)+\sum_{r=1}^\infty\big(S_{i,N,2^r}(t)
 -S_{i,N,2^{r-1}}(t)\big).
 \end{equation}
 In order to complete this programm it remains only to compute limiting 
 covariances as in Section 4 of \cite{KV} taking care also of the slow
 time $s/N$ entering (\ref{3.3}) and (\ref{3.15}).

 \section{Limiting covariances}\label{sec4}
\setcounter{equation}{0}

 In this section we show the existence and compute the limit as $N\to\infty$
 of the expression
 \begin{equation}\label{4.1}
 E\big(G_{i,l}^\ve(s)G_{j,m}^\ve(t)\big)=\ve\int_0^{\tau_i(s/\ve)}
 \int_0^{\tau_j(t/\ve)}b^{l,m}_{ij}(\bar Z(\ve u),\bar Z(\ve v);u,v)dudv.
 \end{equation}
 We start with showing that there exists a constant $C>0$ such that for
 all $t\geq s>0,\, l=1,...,d$, $N\geq 1$ and $i=1,...,\ell$,
 \begin{equation}\label{4.2}
 \sup_{\ve>0}E|G^{1/N}_{i,l}(t)-G^{1/N}_{i,l}(s)|^2\leq C(t-s).
 \end{equation}
 In order to obtain (\ref{4.2}) we note that by (\ref{2.9}) and (\ref{2.10})
 for $t\geq s$,
 \begin{equation}\label{4.3}
 q_i(t)-q_i(s)\geq\al_i(t-s)\,\,\mbox{and}\,\, q_i(t)-q_{i-1}(t)\geq
 \al_{i-1}t\,\,\mbox{when}\,\, i=2,...,k
 \end{equation}
and for any $\gam>0$ there exists $t_\gam$ such that for all $t\geq t_\gam$
and $i=k+1,...,\ell$,
\begin{equation}\label{4.4}
 q_i(t)-q_i(s)\geq (t-s)+\gam^{-1}\,\,\mbox{and}\,\, q_i(t)-q_{i-1}(t)\geq
 t+\gam^{-1}.
 \end{equation}
 Now (\ref{4.2}) follows from (\ref{2.6}), (\ref{3.14}), (\ref{4.1}), 
 (\ref{4.3}) and (\ref{4.4}). Observe, that by (\ref{3.2}) and (\ref{4.1})
 if $\frac 1{N+1}\leq\ve\leq\frac 1N$ then
 \[
 |EG^\ve_{i,l}(s)G^\ve_{j,m}(t)-EG^{1/N}_{i,l}(s)G^{1/N}_{j,m}(t)|\leq
 \frac {4KdC\sqrt \cT}{\sqrt N},
 \]
 and so it suffices to study (\ref{4.1}) as $\ve=\frac 1N$ and $N\to\infty$.
 
 Next, we claim that if $i>j$ and $i>k$ then the limit in (\ref{4.1}) as
 $\frac 1\ve=N\to\infty$ exists and equals zero. Indeed, in this case for 
 any small $\gam >0$ with $\gam T\leq s$,
 \begin{equation}\label{4.5}
 |EG^{1/N}_{i,l}(s)G^{1/N}_{j,m}(t)|\leq I_1+I_2
 \end{equation}
 where by (\ref{4.2}),
 \begin{equation}\label{4.6}
 I_1=|EG^{1/N}_{i,l}(\gam \cT)G^{1/N}_{j,m}(t)|\leq\big(E(G^{1/N}_{i,l}
 (\gam \cT))^2\big)^{1/2}\big(E(G^{1/N}_{j,m}(t))^2\big)^{1/2}\leq C
 \sqrt {\gam \cT t}
 \end{equation}
 and by (\ref{3.14}),
 \begin{eqnarray}\label{4.7}
 &I_2=|E(G^{1/N}_{i,l}(s)-G^{1/N}_{i,l}(\gam \cT))G^{1/N}_{j,m}(t)|=\frac 1N
 \int_{\gam \cT N}^{sN}du\int_0^{\tau_j(tN)}\\
 &b^{l,m}_{ij}(\bar Z(u/N),\bar Z(v/N);u,v)dv\leq\frac CN\int_{\gam \cT N}^{sN}
 du\int_0^{\tau_j(tN)}\rho_{ij}(u,v)dv\nonumber
 \end{eqnarray}
 where
 \begin{equation}\label{4.8}
 \rho_{ij}(u,v)=(\eta_{p,\ka}(n_{ij}(u,v)))^{1-\frac {\vr}{p\te}}
+(\zeta_q(n_{ij}(u,v)))^\del
\end{equation}
with $n_{ij}(s,t)$ defined after (\ref{3.14}). It follows from (\ref{2.6}),
(\ref{2.9}), (\ref{2.10}) and (\ref{4.8}) that for any $\gam>0$ there exists
$N_\gam$ such that whenever $N\geq N_\gam$ and $v\in[0,\cT N]$ (cf. 
Proposition 4.5 in \cite{KV}),
\[
\int_{\gam \cT N}^{sN}\rho_{ij}(u,v)du\leq\gam,
\]
and so $I_2\leq C\cT\gam$. Since $\gam>0$ is arbitrary this together with
(\ref{4.5}) and (\ref{4.6}) yields that for all $l,m=1,...,d$, $i>k$ and
$j<i$,
\begin{equation}\label{4.9}
\lim_{N\to\infty}EG^{1/N}_{i,l}(s)G^{1/N}_{j,m}(t)=0.
\end{equation}

Next, we claim that when $i>k$ then also for all $l,m=1,...,d$,
\begin{equation}\label{4.10}
\lim_{N\to\infty}EG^{1/N}_{i,l}(s)G^{1/N}_{i,m}(t)=0.
\end{equation}
Indeed, by (\ref{3.14}) and (\ref{4.8}) for $t\geq s$,
\begin{equation}\label{4.11}
|EG^{1/N}_{i,l}(s)G^{1/N}_{i,m}(t)|\leq\frac 1N\int_0^{sN}du\int_0^{tN}
\rho_{ii}(u,v)dv=I_3+I_4
\end{equation}
where 
\[
I_3=\frac 2N\int_0^{sN}du\int_u^{sN}\rho_{ii}(u,v)dv\,\,\mbox{and}\,\,
I_4=\frac 1N\int_0^{sN}du\int_{sN}^{tN}\rho_{ii}(u,v)dv.
\]
Now
\begin{eqnarray}\label{4.12}
&I_3=\frac 2N\int_0^{sN}du\int_u^{u+\gam}\rho_{ii}(u,v)dv+\frac 2N
\int_0^{\gam N}du\int_{u+\gam}^{sN}\rho_{ii}(u,v)dv\\
&+\frac 2N\int_{\gam N}^{sN}du\int_{u+\gam}^{sN}\rho_{ii}(u,v)dv\leq C(s\gam
+\gam+s\be_\gam(\gam N))\nonumber
\end{eqnarray}
for some $C>0$ where by (\ref{2.6}) and (\ref{2.10}) for any $\gam>0$,
\begin{equation}\label{4.13}
\be_\gam(M)=\sup_{u\geq M}\int^\infty_{u+\gam}\rho_{ii}(u,v)dv<\infty\,
\,\mbox{and}\,\, \lim_{M\to\infty}\be_\gam(M)=0.
\end{equation}

Next,
\begin{eqnarray}\label{4.14}
&I_4=\frac 1N\int_0^{sN}du\int_{sN}^{sN+\gam}\rho_{ii}(u,v)dv\\
&+\frac 1N\int_0^{sN}du\int_{sN+\gam}^{tN}\rho_{ii}(u,v)dv\leq Cs\gam+
Cs\be_s(N).\nonumber
\end{eqnarray}
Finally, (\ref{4.10}) follows from (\ref{4.11})--(\ref{4.14}) letting, first,
$N\to\infty$ and then $\gam\to 0$.

In order to compute the limit as $\frac 1\ve=N\to\infty$ of (\ref{4.1}) 
for $i,j=1,2,...,
k$ we recall an argument of Lemma 3.1 from \cite{Kh1} which yields that if
uniformly in $\sig\geq 0$ and $x,y$ from a compact set the limit
\begin{equation}\label{4.15}
\lim_{N\to\infty}\frac 1N\int_{\sig/\al_i}^{(\sig+sN)/\al_i}du
\int_{\sig/\al_j}^{(\sig+sN)/\al_j}b_{ij}^{lm}(x,y;u,v)dudv=sD^{l,m}_{ij}(x,y)
\end{equation}
exists and has the form of the right hand side with a continuous
$D^{l,m}_{i,j}$ then the limit (\ref{4.1}) exists, as well, and it has the 
form
\begin{equation}\label{4.16}
\lim_{N\to\infty}E(G^{1/N}_{i,l}(s)G^{1/N}_{j,m}(t))=\int_0^{\min(s,t)}
D^{l,m}_{ij}(\bar Z(u),\bar Z(u))du.
\end{equation}

Namely, set $M=M(N)=[N^{2/3}]$ and let $s_\iota=\frac {\iota s}M,\, 
\iota=0,1,...,M-1$. Assume also that $s\leq t$. Let
\[
A_N=\cup_{\iota=0}^{M-1}A_{N,\iota}\,\,\mbox{with}\,\, A_{N,\iota}=\{ (u,v):\,
s_\iota N\leq u,v<(s_\iota+\frac sM)N\}
\]
and $B_N=\{(u,v):\, 0\leq u\leq sN,\, 0\leq v\leq tN\}\setminus A_N$. Then
\begin{equation}\label{4.17}
EG^{1/N}_{i,l}(s)G^{1/N}_{j,m}(t)=I_5+I_6
\end{equation}
where
\[
I_5=\frac 1{Nij}\int_{B_N}b^{lm}_{ij}(\bar Z(u/N),\bar Z(v/N),u/i,v/j)dudv
\]
and
\[
I_6=\frac 1{Nij}\int_{A_N}b^{lm}_{ij}(\bar Z(u/N),\bar Z(v/N),u/i,v/j)dudv.
\]
Now, by (\ref{3.14}) and (\ref{4.8}),
\begin{eqnarray}\label{4.18}
&|I_5|\leq\frac C{Nij}\big(\sum_{\iota=0}^{M-1}\int_0^{s_\iota/\ve}
\int_{s_\iota/\ve}^{(s_\iota+\frac sM)/\ve}(\rho_{ij}(u/\al_i,v/\al_j)\\
&+\rho_{ji}(u/\al_j,v/\al_i))dudv+\int_0^{s_\iota/\ve}
\int_{s_\iota/\ve}^{t_\iota/\ve}\rho_{ij}(u/\al_i,v/\al_j)dudv\big).\nonumber
\end{eqnarray}

Observe that by the definition of $n_{ij}(u,v)$ after (\ref{3.14}) we can
write for $i,j=1,...,k$,
\begin{equation}\label{4.19}
\rho_{ij}(u/\al_i,v/\al_j)=\zeta(|u-v|)
\end{equation}
where $\zeta\geq 0$ satisfies $\int_0^\infty w\zeta(w)dw<\infty$. Integrating
by parts we obtain for any $V\geq U\geq 0$,
\begin{equation}\label{4.20}
\int_0^Udu\int_U^V\zeta(v-u)dv\leq\int_0^Udu\int_{U-u}^\infty\zeta(w)dw=
\int_0^Ur\zeta(r)dr\leq\int_0^\infty r\zeta(r)dr.
\end{equation}
This together with (\ref{2.6}), (\ref{4.18}) and (\ref{4.19}) gives by the
choice of $M=M(N)$ that 
\begin{equation}\label{4.21}
|I_5|\leq \tilde C\frac M{N\al_i\al_j}\to 0\,\,\mbox{as}\,\, N\to\infty
\end{equation}
for some $\tilde C>0$ independent of $M$ and $N$.

Next,
\begin{equation}\label{4.22}
I_6=\frac 1{M\al_i\al_j}\sum_{\iota=0}^{M-1}J_{M,N}(\iota)+I_7
\end{equation}
where
\[
J_{M,N}(\iota)=\frac MN\int_{s_\iota\leq u,v<(s_\iota+\frac sM)N}b_{ij}^{lm}
(\bar Z(s_\iota),\bar Z(s_\iota);\frac u{\al_i},\frac v{\al_j})dudv
\]
and by (\ref{2.7}) and the choice of $M=M(N)$,
\begin{equation}\label{4.23}
|I_7|\leq Cs^3NM^{-2}\to 0\,\,\mbox{as}\,\, N\to\infty
\end{equation}
where $C>0$ does not depend on $s,N$ and $M$. By (\ref{4.15}) we obtain that
\begin{equation}\label{4.24}
|J_{M,N}(\iota)-s\al_i\al_jD_{ij}^{l,m}(\bar Z(s_\iota),\bar Z(s_\iota))|\to
0\,\,\mbox{as}\,\, N\to\infty,
\end{equation}
and so
\begin{equation}\label{4.25}
|I_6-\int_0^sD_{ij}^{l,m}(\bar Z(u),\bar Z(u))du|\to
0\,\,\mbox{as}\,\, N\to\infty
\end{equation}
completing the proof of (\ref{4.16}).

In order to describe $D^{l,m}_{ij}(x,y),\, i,j\leq k$ consider all indices
$1\leq i'_1<i'_2<...<i'_{\iota_{ij}}=i$ and $1\leq j'_1<j'_2<...
<j'_{\iota_{ij}}=j$ such that there exist $0<\rho_1<...<\rho_{\iota_{ij}}=1$
satisfying $\al_{i'_l}\rho_l,\,\al_{j'_l}\rho_l\in\{\al_1,...,\al_k\}$ for
all $l=1,...,\iota_{ij}$.
Define
\begin{eqnarray}\label{4.26}
&a^{l,m}_{ij}(x,y;s_1,...,s_{\iota_{ij}})=\int B^{(l)}_i(x,\xi_1,...,\xi_i)
B_j(y,\tilde\xi_1,...,\tilde\xi_j)\\
&\prod_{\be=1}^{\iota_{ij}}d\mu_{s_\be}(\xi_{i'_\be},\tilde\xi_{j'_\be})
\prod_{i_\gam\not\in\{i'_1,...,i'_{\iota_{ij}}\},1\leq i_\gam<i}
d\mu(\xi_{i_\gam})
\prod_{j_\zeta\not\in\{j'_1,...,j'_{\iota_{ij}}\},1\leq j_\zeta<j}
d\mu(\xi_{j_\zeta}).\nonumber
\end{eqnarray}
Then in the same way as in the proof of Lemma 4.4 from \cite{KV} (see also
Section 7 there) we obtain relying on (\ref{2.6}), (\ref{3.10}) and 
(\ref{3.14}) that
\begin{equation}\label{4.27}
\lim_{N\to\infty,\,\al_iNu_N-\al_jNv_N=w}b_{ij}^{l,m}(x,y;Nu_N,Nv_N)=
a_{ij}^{l,m}(x,y;\rho_1w,\rho_2w,...,\rho_{\iota_{ij}}w).
\end{equation}
This is the only place where we need Assumption \ref{ass2.1} for 
$\eta_{p,\ka,s}$ with $s>0$.
It follows similarly to Section 7 of \cite{KV} that the limit (\ref{4.15})
 exists and it can be written in the form
 \begin{equation}\label{4.28}
 D^{l,m}_{ij}(x,y)=\frac 1{\al_i\al_j}\int_{-\infty}^\infty
 a_{ij}^{l,m}(x,y;\rho_1w,\rho_2w,...,\rho_{\iota_{ij}}w)dw.
 \end{equation}
 
 Collecting the results of Sections \ref{sec3} and \ref{sec4} together we 
 conclude that each $G_i^\ve,\, i=1,...,k$ converges weakly as $\ve\to 0$
  to the corresponding Gaussian process $G^0_i$ having independent increments
  while the process $G^\ve_i,\, i>k$ converge weakly as $\ve\to 0$ to zero
  (in the continuous time case we are dealing with now). Moreover, the 
  processes $G^\ve$ converge weakly as $\ve\to 0$ to a Gaussian process
  $G^0$ (with not necessarily independent increments as an example in 
  \cite{KV} shows) having the representation
  \begin{equation}\label{4.29}
  G^0(t)=\sum_{i=1}^kG^\ve_i(it).
  \end{equation}
  Furthermore, the covariances of different components $G^0_i(s)=(G^{0,1}_i,
  ...,G^{0,d}_i(s))$ of this processes are described in view of the above by
  \begin{equation}\label{4.30}
  EG^{0,l}_i(s)G^{0,m}_j(t)=\int_0^{\min(s,t)}D^{l,m}_{ij}(\bar Z(u),
  \bar Z(u))du,
  \end{equation}
  and so by (\ref{4.29}),
  \begin{equation}\label{4.31}
  EG^{0,l}(s)G^{0,m}(t)=\int_0^{\min(s,t)}A^{l,m}(u)du
  \end{equation}
  where
  \[
  A^{l,m}(u)=\sum_{1\leq i,j\leq k}D_{ij}^{l,m}(\bar Z(iu),\bar Z(ju)).
  \]

\section{Gaussian approximation of the slow motion and discrete time case}
\label{sec5}\setcounter{equation}{0}

In order to complete the proof of Theorem \ref{thm2.2} we proceed similarly 
to \cite{Kh1}. First, we consider the process $H^\ve(t)$ which solves the
linear equation
\begin{equation}\label{5.1}
H^\ve(t)=G^\ve(t)+\int_0^t\nabla\bar B(\bar Z(s))H^\ve(s)ds.
\end{equation}
By (\ref{2.7}), for some $C>0$ independent of $t$ and $\ve$,
\[
|H^\ve(t)|\leq |G^\ve(t)|+C\int_0^t|H^\ve(s)|ds.
\]
Then
\[
\big\vert |H^\ve(t)|-|G^\ve(t)|\big\vert\leq C\int_0^t|G^\ve(s)|ds+
C\int_0^t\big\vert |H^\ve(s)|-|G^\ve(s)|\big\vert ds
\]
and by Gronwall's inequality
\begin{equation}\label{5.2}
|H^\ve(t)|\leq |G^\ve(t)|+Ce^{Ct}\int_0^t|G^\ve(s)|ds.
\end{equation}
It follows from Section \ref{sec3} that the family of processes $\{ G^\ve(t),
\, t\in[0,\cT]\}$ is tight which together with (\ref{5.2}) implies that the
family of processes $\{ H^\ve(t),\, t\in[0,\cT]\}$, as well, as the family
of pairs $V^\ve=\{ G^\ve,\, H^\ve\}$ are tight. 

It follows that any weak limit $V^0=\{ G^0,\, H^0\}$ of $V^\ve$ as 
$\ve\to 0$ must satisfy the equation
\begin{equation}\label{5.3}
H^0(t)=G^0(t)+\int_0^t\nabla\bar B(\bar Z(s))H^0(s)ds
\end{equation}
which has a unique solution. Moreover, its solution $H^0$ is a Gaussian
process. Indeed, the equation (\ref{5.3}) can be
solved by successive approximations starting from $G^0$ so that on each 
step we will get a Gaussian process (in view of linearity) and the 
limiting process will be Gaussian, as well. Moreover, $H^0$ depends linearly
 on $G^0$ having an integral representation of the form
\begin{equation}\label{5.4}
H^0(t)=G^0(t)+\int_0^tK(t,s)G^0(s)ds
\end{equation}
with a differentiable kernel $K$ (Green's function). The latter follows
considering an operator $A$ given by 
\[
Af(t)=\int_0^t\nabla\bar B(\bar Z(s))f(s)ds
\]
which has the supremum norm less than 1 if $t\in[0,\Del]$ for $\Del$ small
enough, and so we can write 
\[
H^0=(I-A)^{-1}G^0=G^0+\sum_{n=1}^\infty A^nG^0.
\]
In view of the form of the integral operator $A$ above this representation
yields (\ref{5.4}) on the interval $[0,\Del]$ and then employing the same
argument successively to time itervals $[\Del,2\Del],\,[2\Del,3\Del],...$
we extend the representation (\ref{5.4}) for any $t$.

Observe that
\begin{eqnarray}\label{5.5}
&Q^\ve(t)=\ve^{-1/2}\int_0^t\big( B(Z^\ve_x(s),\xi(q_1(s/\ve)),...,
\xi(q_\ell(s/\ve)))-\bar B(\bar Z_x(s))\big)ds\\
&=G^\ve(t)+\int_0^t\nabla_x B(Z^\ve_x(s),\xi(q_1(s/\ve)),...,
\xi(q_\ell(s/\ve)))Q^\ve(s)ds+\int_0^tJ^\ve_1(s)ds\nonumber
\end{eqnarray}
where
\begin{eqnarray*}
&J^\ve_1(s)=\ve^{-1/2}\big( B(\bar Z_x(s)+\sqrt\ve Q^\ve(s),\xi(q_1(s/\ve))
,...,\xi(q_\ell(s/\ve)))-B(\bar Z_x(s),\\
&\xi(q_1(s/\ve)),...,\xi(q_\ell(s/\ve)))-\nabla_xB(\bar Z_x(s),
\xi(q_1(s/\ve)),...,\xi(q_\ell(s/\ve)))\sqrt\ve Q^\ve(s)\big).
\end{eqnarray*}
If $H^\ve$ solves (\ref{5.1}) then $U^\ve(t)=Q^\ve(t)-H^\ve(t)$ satisfies
by (\ref{5.4}) the equation
\begin{equation}\label{5.6}
U^\ve(t)-\int_0^t\nabla_xB(\bar Z_x(s),\xi(q_1(s/\ve)),...,\xi(q_\ell(s/\ve)))
U^\ve(s)ds=\int_0^t(J^\ve_1(s)+J^\ve_2(s))ds
\end{equation}
where
\[
J^\ve_2(s)=\big(\nabla_xB(\bar Z_x(s),\xi(q_1(s/\ve)),...,\xi(q_\ell(s/\ve)))
-\nabla_x\bar B(\bar Z_x(s))\big)H^\ve(s).
\]
By Gronwall's inequality we obtain that
\begin{equation}\label{5.7}
|U^\ve(t)|\leq Cte^{ct}\int_0^t|J^\ve_1(s)+J^\ve_2(s)|ds
\end{equation}
for some $C>0$ independent of $\ve$ and $t\in[0,\cT]$.

Thus, in order to prove that $Q^\ve$ converges weakly as $\ve\to 0$ to
a Gaussian process $Q^0$ solving (\ref{2.14}) it suffices to show that
$\int_0^tJ^\ve_1(s)ds$ and $\int_0^tJ^\ve_2(s)ds$ converge to zero in
probability as $\ve\to 0$. By (\ref{2.7}),
\[
|Z^\ve_x(t)-Y^\ve_x(t)|\leq C\int_0^t|Z^\ve_x(s)-\bar Z_x(s)|ds=C\sqrt\ve
\int_0^t|Q^\ve_x(s)|ds
\]
with $C=Kd$, and so
\[
|Q^\ve_x(t)|\leq |G^\ve(t)|+C\int_0^t|Q^\ve_x(s)|ds.
\]
Hence, in the same way as in (\ref{5.2}),
\begin{equation}\label{5.8}
|Q^\ve_x(t)|\leq |G^\ve(t)|+Ce^{Ct}\int_0^t|G^\ve(s)|ds.
\end{equation}
By (\ref{2.7}) and the Taylor formula with a reminder we conclude that
\begin{equation}\label{5.9}
|J^\ve_1(s)|\leq C\sqrt\ve |Q^\ve(s)|^2
\end{equation}
which together with (\ref{4.2}) yields that $E|J^\ve_1(s)|\to 0$ as 
$\ve\to 0$.

The proof of convergence to zero in probability of $\int_0^tJ^\ve_2(s)ds$ 
as $\ve\to 0$ is based on the integral representation (\ref{5.4}). Set
\[
\Phi(x,\xi_1,...,\xi_\ell)=B(x,\xi_1,...,\xi_\ell)-\bar B(x)
\]
and
\[
\Psi(x,\xi_1,...,\xi_\ell)=\nabla_xB(x,\xi_1,...,\xi_\ell)-\nabla_x\bar B(x).
\]
Relying on the representation (\ref{5.4}) we obtain that
\begin{equation}\label{5.10}
\big\vert E\int_0^tJ_2^\ve(s)ds\big\vert\leq |J^\ve_3(t)|+|J^\ve_4(t)|
\end{equation}
where
\begin{eqnarray}\label{5.11}
&J^\ve_3(t)=\ve^{3/2}\int_0^{t/\ve}ds\int_0^sduE\big(\Psi(\bar Z_x(\ve s),
\xi(q_1(s)),...,\xi(q_\ell(s)))\\
&\times \Phi(\bar Z_x(\ve u),\xi(q_1(u)),...,\xi(q_\ell(u)))\big)\nonumber
\end{eqnarray}
and
\begin{eqnarray}\label{5.12}
&J^\ve_4(t)=\ve^{3/2}\int_0^{t/\ve}ds\int_0^{\ve s}du\int_0^{u/\ve}dv
K(\ve s,\ve v)\\
&\times E\big(\Psi(\bar Z_x(\ve s),\xi(q_1(s)),...,\xi(q_\ell(s)))
\Phi(\bar Z_x(\ve u),\xi(q_1(u)),...,\xi(q_\ell(u)))\big).\nonumber
\end{eqnarray}
Estimating the expectations in (\ref{5.11}) and (\ref{5.12}) via 
(\ref{3.10}) similarly to (\ref{3.14}) we obtain that both $J_3^\ve(t)$
and $J^\ve_4(t)$ are of order $\sqrt\ve$, and so the left hand side of
(\ref{5.10}) is of order $\sqrt\ve$, as well. For more details of a
similar argument we refer the reader to \cite{Kh1}. This completes the
proof of Theorem \ref{2.2} concerning the continuous time case.

In the discrete time case the proofs are similar but slightly simpler.
Namely, set
\begin{equation}\label{5.13}
R^{1/N}_i(t)=N^{-1/2}\sum_{n=0}^{[Nt/i]}B_i(\bar Z(nt/N),\xi(q_1(n)),
...,\xi(q_i(n)))
\end{equation}
where $B_i$'s are the same as in (\ref{2.16})--(\ref{2.18}). Then for 
all $N\geq 1$,
\begin{equation}\label{5.14}
|G^{1/N}_i(t)-R^{1/N}_i(t)|\leq CN^{-1/2}
\end{equation}
for some $C>0$ independent of $N$. The asymptotical behavior of $R^{1/N}$
as $N\to\infty$ can be studied in the same way as in \cite{KV} taking 
into account that we have here slightly different mixing conditions, and so 
the corresponding estimates should be done as above via 
(\ref{3.10})--(\ref{3.14}). The main difference of the discrete vis-\' a-vis
continuous time case is that now each $G^{1/N}_i(t),\, i=k+1,...,\ell$ converges
weakly as $N\to\infty$ to a nondegenerate Gaussian process $G^0_i(t)$ having
the covariances
\begin{equation}\label{5.15}
E\big(G^0_i(t)G^0_i(s)\big)=\int_0^{\min(s,t)}du\int\big(B_i(\bar Z(u),
\xi_1,...,\xi_i)\big)^2d\mu(\xi_1)...d\mu(\xi_i)
\end{equation}
which is proved combining arguments of Proposition 4.5 in \cite{KV} and
of Section \ref{sec4} above. The computation of other limiting covariances
proceeds in the same way as in the continuous time case.
It follows that in the discrete time case the processes $G^\ve$ converge
weakly as $\ve\to 0$ to a Gaussian process $G^0$ having the representation
\[
G^0(t)=\sum_{i=1}^k G^0_i(it)+\sum_{i=k+1}^\ell G^0_i(t)
\]
where each process $G^0_i,\, i>k$ is independent of each $G^0_j$ with
$j\ne i$ while the processes $G^0_i,\, i\leq k$ are correlated with
covariances described at the end of Section \ref{sec4} taken with 
$\al_i=i,\, i=1,...,k$. The argument concerning the convergence of 
processes $Q^\ve$ to $Q^0$ solving (\ref{2.14}) remains the same as
in the continuous time case.

\section{Some dynamical systems applications}\label{sec6}\setcounter{equation}
{0}

We start with recalling the setup from \cite{Do1} and \cite{Do2}. A 
$C^2$-diffeomorphism $F$ of a compact Riemannian manifold $\Om$ is called
partially hyperbolic if there is a $F$-invariant splitting $E^u\oplus E^c
\oplus E^s$ of the tangent bundle of $\Om$ with $E^u\ne 0$ and constants
$\la_1\leq\la_2<\la_3\leq\la_4<\la_5\leq\la_6$, $\la_2<1,\,\la_5>1$ such 
that $\|dF(v)\|/\|v\|$ is between $\la_1$ and $\la_2$ on $E^s$, between 
$\la_3$ and $\la_4$ on $E^c$ and between $\la_5$ and $\la_6$ on $E^u$.
Denote by $W^u$ the foliation tangent to $E^u$ and call $S$ a u-set if $S$
belongs to a single leaf of $W^u$. $F$-invariant probability measures
which are absolutely continuous with respect to the volume on leafs $W^u$
are called u-Gibbs measures. It is assumed that $F$ has a unique u-Gibbs
measure $\SRB$ which is called the Sinai-Ruelle-Bowen (SRB) measure.

An important role in the construction is played by Markov families which
 are collections $\cS$ of u-sets which cover $\Om$ and have certain 
 regularity properties (see \cite{Do1} and \cite{Do2}) but we formulate
  here only their "Markov property" saying that for any $S\in\cS$ there 
  are $S_i\in\cS$ such that $FS=\cup_iS_i$. Now let $\cS$ be a Markov
  family. Following \cite{Do1} and \cite{Do2} we construct on each $S\in\cS$
   an increasing sequence of $\sig$-algebras $\cF_n$ in 
   the following recursive way. Let $\cF^S_0=\{\emptyset,S\}$. Suppose that
   $\cF_n^S$ is generated by $\{ S_{j,n}\}$ with $F^nS_{j,n}\in\cS$. By
   the "Markov property" we can decompose $F^{n+1}S_{j,n}=\cup_lS_{jl,n}$
   and now let $\cF^S_{n+1}$ be generated by $F^{-n-1}S_{jl,n}$.
   
   Next, for each $x_1$ and $x_2$ in a u-set $S$ put
   \[
   \rho(x_1,x_2)=\prod_{j=0}^\infty\frac {\mbox{det}(dF^{-1}|E^u)(F^{-1}x_1)}
   {\mbox{det}(dF^{-1}|E^u)(F^{-1}x_2)}.
   \]
   Fix $x_0\in S$ and let $\rho_S(x)=\rho(x,x_0)(\int_S\rho(x,x_0)dx)^{-1}$.
   For a Markov family $\cS$ and nonnegative constants $R,\al$ denote by
   $E_1(\cS,R,\al)$ the set of probability measures $\sig$ defined for each
    continuous function $g\in C(\Om)$ by
   \begin{equation}\label{6.1}
   \sig(g)=\int_Sg(x)e^{G(x)}\rho_S(x)dx
   \end{equation}
   where $S\in\cS$ and $G$ is H\" older continuous with the exponent $\al$
    and the constant $R$. Denote also by $E=E(\cS,R,\al)$ the closure of the
    convex hull of $E_1(\cS,R,\al)$. The decay of correlations 
    is measured in \cite{Do1} and \cite{Do2} via a sequence $a(n)\to 0$ as
    $n\to\infty$ such that for any $\sig\in E$ and each H\" older continuous 
    $g$ on $\Om$,
    \begin{equation}\label{6.2}
    |\sig(g\circ F^n)-\SRB(g)|\leq a(n)\| g\|
    \end{equation}
    where $\|\cdot\|$ is a H\" older norm. An argument from Section 5
    of \cite{DL} compares the coefficient $a(n)$ above with the more familiar
    rate of decay of correlations $|\SRB(f\cdot(g\circ F^n))-\SRB(f)\SRB(g)|$
    and it follows from there that the latter decays superpolynomially if
    and only if $a(n)$ decays superpolynomially. According to \cite{Do0}
    such decay of correlations holds true for $C^2$ Anosov flows with
    jointly nonintegrable stable and unstable foliations and for their
    time-one maps. By \cite{FMT}
    this remains true for an open dense set of $C^2$ Axion A flows as well,
    as for their time-one maps. For other partially hyperbolic dynamical
    systems with fast decay of correlations see \cite{Do1}, \cite{Do2},
    \cite{FMT} and references there.
    
    In order to estimate $\eta_{p,\ka,s}(n)$ from (\ref{2.3}) we write in
    the same way as in Lemma 4 from \cite{Do1} that on each element $S$
    in $\cF_{[t]}$,
    \begin{equation}\label{6.3}
    A_{n,s,t}=E\big(g(f\circ F^{n+t},f\circ F^{n+t+s})|\cF_{[t]}\big)=
    \int_S\rho_S(y)g_{s,t}(F^ny)dy
    \end{equation}
    where the expectation is with respect to $\sig$ on $S$ and $g_{s,t}(z)=
    g(f(F^{t-[t]}z),f(F^{t-[t]+s}z))$. If $f$ and $g$ are H\" older 
    continuous then $g_{s,t}$
    is H\" older continuous for fixed $s$ and $t$ and it is uniformly in $t$
    H\" older continuous when $s=0$. Thus, by (\ref{6.2}) we have that
    $|A_{n,s,t}-EA_{n,s,t}|$ decays in $n$ with the speed of at least $a(n)$
    and this decay is uniform in $t$ if $s=0$. Hence, if $a(n)$ decays
    superpolynomially then (\ref{2.6}) holds true. This yields Theorem
    \ref{thm2.2} for $\xi(t)=\xi(t,z)=g(F^tz)$ on a probability space
    $(S,\sig)$ for $\sig\in E$ and an element $S$ of a Markov family while
    $g$ is a H\" older continuous function. We observe that the measure
    $\sig$ here plays the role of the probability $Pr$ in the setup of
    Section \ref{sec2} while $\SRB$ plays the role of $P$ there.

    \section{Concluding remarks: fully coupled averaging}\label{sec7}
\setcounter{equation}{0}

In the nonconventional framework as discussed in this paper even the setup
of fully coupled averaging, i.e. when the fast motion depends on the slow
one, is not quite clear.
On the first sight we may want to deal with the equations
\begin{eqnarray}\label{7.1}
&X^\ve(n+1)=X^\ve(n)+\ve B(X^\ve(n),\xi(n),\xi(2n),...,\xi(\ell n)),\\
&\xi(n+1)=F_{X^\ve(n)}(\xi(n))\nonumber
\end{eqnarray}
in the discrete time case and
\begin{equation}\label{7.2}
\frac {dX^\ve(t)}{dt}=\ve B(X^\ve(t),\xi(t),\xi(2t),...,\xi(\ell t)),\quad
\frac {d\xi(t)}{dt}=b(X^\ve(t),\xi(t))
\end{equation}
in the continuous time case. The problem is that $\xi(kn)$ or $\xi(kt)$ are not
yet defined for $k>1$ at time $n$ or $t$ so we cannot insert them into the 
first equation in (\ref{7.1}) or (\ref{7.2}) respectively, and so these 
equations do not define properly $X^\ve$ and $\xi$.

A reasonable modification of this setup is to consider
\begin{eqnarray}\label{7.3}
&X^\ve(n+1)=X^\ve(n)+\ve B(X^\ve(n),\eta_1(n),\eta_2(n),...,\eta_\ell(n)),\\
&\eta_i^\ve(n+1)=F^i_{X^\ve(n)}(\eta_i^\ve (n)),\, i=1,...,\ell\nonumber
\end{eqnarray}
in the discrete time case and
\begin{eqnarray}\label{7.4}
&\frac {dX^\ve(t)}{dt}=\ve B(X^\ve(t),\eta_1(t),\eta_2(t),...,\eta_{\ell}(t)),\\
&\frac {d\eta_i^\ve(t)}{dt}=ib(X^\ve(t),\eta^\ve_i(t)),\, i=1,2,...,\ell
\nonumber\end{eqnarray}
in the continuous time case. We consider (\ref{7.3}) and (\ref{7.4}) as sets 
of $\ell+1$ equations but require that $\eta^\ve_1(0)=\eta^\ve_2(0)=\cdots=
\eta^\ve_\ell(0)$. This approach seems to be reasonable if we consider 
(\ref{7.3}) and (\ref{7.4}) as perturbations of equations with constants of 
motion 
\begin{equation}\label{7.5}
 \eta^{(x)}(n+1)=F_x(\eta^{(x)}(n))\,\,\mbox{and}\,\,\frac {d\eta^{(x)}(t)}{dt}
 =B(x,\eta^{(x)}(t)),
\end{equation}
i.e. when $x$ variable remains fixed in unperturbed equations but start
moving slowly in perturbed ones.
Then $\eta^{(x)}(i(n+1))=F^i_x(\eta^{(x)}(in))$ and $d\eta^{(x)} (it)/dt=
iB(x,\eta^{(x)}(it))$.

As it is well known in the fully coupled setup the averaging principle not 
always holds true and when it takes place then usually only in the sense of 
convergence in average or in measure. In the nonconventional situation the 
problem is even more complicated. Consider, for instance,
\begin{eqnarray}\label{7.6}
&\frac {d\al^\ve_{\al,\vf}(t)}{dt}=\ve B(\al^\ve_{\al,\vf}(t),\vf_{1,\al}^\ve(t)
,...,\vf_{\ell,\al}^\ve(t)),\\
&\frac {d\vf_{i,\al,\vf}^\ve(t)}{dt}=i\al^\ve_{i,\al,\vf}(t),\quad \al^\ve_{\al,
\vf}(0)=\al,\,\vf_{1,\al,\vf}^\ve(0)=\cdots=\vf_{\ell,\al,\vf}^\ve(0)=\vf
\nonumber\end{eqnarray}
where $\vf$ denotes a point on an $n$-dimensional torus $\bbT^n$ and $\al$ 
denotes a constant $n$-vector (constant vector field on $\bbT^n$). 
 Then $\vf^\ve_{i,\al,\vf}=i\vf^\ve_{1,\al,\vf}-
(i-1)\vf$. Set $\tilde B(\psi,\vf)=B(\al,\psi,\psi-\vf...,\psi-(\ell-1)\vf)$. 
Then the right
hand side of (\ref{7.6}) can be replaced by $\ve \tilde B(\al^\ve_{\al,\vf}(t),
\vf_{1,\al}^\ve(t),\vf)$.
If $\bar B(\al)=\int\tilde B(\al,\vf_1,\vf)d\vf_1d\vf$ and
$\frac {d\bar\al_{\al}(t)}{dt}=\bar B(\bar\al_{\al}(t)),\quad \bar\al_{\al}(0)
=\al$ then employing the technique from the proof of Theorem 2.1 in
\cite{Ki2} it is not difficult to see that for any compact $K$, 
\begin{equation}\label{7.7}
\int_K\sup_{0\leq t\leq\cT/\ve}|\al^\ve_{\al,\vf}(t)-\bar \al_{\al}(\ve t)|
d\al d\vf\to 0\,\,\,\mbox{as}\,\,\, \ve\to 0.
\end{equation}

\bibliography{matz_nonarticles,matz_articles}
\bibliographystyle{alpha}

\end{document}